\documentclass[12pt,reqno]{article}
\usepackage{fullpage}
\usepackage{enumerate}
\usepackage{booktabs}
\usepackage{amsthm,amsmath,amssymb}
\usepackage{tikz}
\usepackage{subfigure}
\usepackage{graphicx}
\usepackage[export]{adjustbox}
\definecolor{webgreen}{rgb}{0,.5,0}
\definecolor{webbrown}{rgb}{.6,0,0}
\usepackage[colorlinks=true,
linkcolor=webgreen,
filecolor=webbrown,
citecolor=webgreen]{hyperref}

\newcommand{\seqnum}[1]{\href{http://oeis.org/#1}{\underline{#1}}}
\newcommand{\Figs}[1]{\hyperref[#1]{Figure~\ref*{#1}}}
\newcommand{\Tabs}[1]{\hyperref[#1]{Table~\ref*{#1}}}
\newcommand{\Secs}[1]{\hyperref[#1]{section~\ref*{#1}}}
\everymath{\displaystyle}

\title{\bf Note on sequences A123192, A137396 and A300453}
\author{Franck Ramaharo\\
\small D\'epartement de Math\'ematiques et Informatique\\[-0.8ex]
\small Universit\'e d'Antananarivo\\[-0.8ex] 
\small 101 Antananarivo, Madagascar\\
\small\href{mailto:franck.ramaharo@gmail.com}{\tt franck.ramaharo@gmail.com}\\
}

\date{\small\today\\}

\begin{document}

\maketitle
\begin{abstract}
 We give the connection between three polynomials that generate triangles in \textit{The On-Line Encyclopedia of Integer Sequences} (\seqnum{A123192}, \seqnum{A137396} and \seqnum{A300453}). We show that they are related with the bracket polynomial for the $(2,n)$-torus knot.
 
\end{abstract}

\section{Introduction}
\begin{figure}[!ht]
    \centering
    \hspace*{\fill}
    \subfigure[]{\includegraphics[width=.215\linewidth]{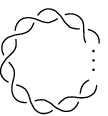}\label{subfig:2ntorusdiagram}}
    \hspace*{\fill}
    \subfigure[]{\includegraphics[width=0.215\linewidth]{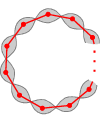}\label{subfig:ncyclegraph}}\hfill%
    \subfigure[]{\includegraphics[width=0.215\linewidth]{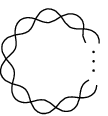}\label{subfig:torusshadow}}%
    \hspace*{\fill}
    \caption{\subref{subfig:2ntorusdiagram} $(2,n)$-torus knot diagram, \subref{subfig:ncyclegraph} $n$-cycle graph and its ``medial graph'',  \subref{subfig:ncyclegraph} $(2,n)$-torus knot shadow diagram.}
    \label{fig:link}
\end{figure}

Let $K_n$ denote the $ (2,n) $-torus knot diagram (\Figs{fig:link}~\subref{subfig:2ntorusdiagram}). The corresponding bracket polynomial is given by the formula \cite{Lafferty}
\begin{equation}\label{eq:bptn}
    K_n(A,B,d) = \dfrac{\left(A + Bd\right)^{n} + \left(d^2 - 1\right)A^n}{d}.
\end{equation}
The following triangles are in \textit{The On-Line Encyclopedia of Integer Sequences}  \cite{Sloane}, and consist of the coefficients in the expansion of $gK_n(A,B,d)$ for some values of $ A,B,d $ and $g$  (see \Tabs{tab:A123192}, \Tabs{tab:A137396} and \Tabs{tab:A300453}).

\begin{itemize}
    \item Row $n$ in \seqnum{A123192} is generated by
    
            \begin{equation}\label{eq:brack}
                x^{|3n-2|}K_n\left(x,x^{-1},-x^{-2}-x^{2}\right) = \begin{cases}
                    \dfrac{\left(x^8 + x^4 + 1\right)x^{4n} + (-1)^nx^4}{x^8 + x^4} & \textit{if}\ n\geq 1;\\
                    -x^4 - 1 & \textit{if}\ n = 0. 
                \end{cases}
            \end{equation}
      
    The interpretation for the choice of $A,B,d$ and $g$ is given in \Secs{sec:interpretation} .
           \begin{table}[!ht]
              \centering
                $\begin{array}{c|rrrrrrrrrrrrrrrrrrrrr}
                n\ \backslash\ k & 0 & 1 & 2 & 3 & 4 & 5 & 6 & 7 & 8 & 9 & 10 & 11 & 12 & 13 & 14 & 15 & 16\\
                \hline
                0 & -1 & 0 & 0 &  0 & -1 \\
                1 &  0 & 0 & 0 &  0 & -1 \\
                2 & -1 & 0 & 0 &  0 &  0 & 0 & 0 & 0 & -1 \\ 
                3 &  0 & 0 & 0 & -1 &  0 & 0 & 0 & 0 &  0 & 0 & 0 & -1 \\
                4 & -1 & 0 & 0 &  0 &  1 & 0 & 0 & 0 & -1 & 0 & 0 &  0 & 0 & 0 & 0 & 0 & -1 
              \end{array}$
              \caption{First $5$ rows in \seqnum{A123192}.}
              \label{tab:A123192}
          \end{table}
    
    \item Row $n$ in \seqnum{A137396} is generated by the chromatic polynomial of the $ n $-cycle graph.
     
        \begin{equation}\label{eq:chrom}
            x^{\frac{1}{2}}K_n\left(-1, x^{\frac{1}{2}}, x^{\frac{1}{2}}  \right) = (x - 1)^n + (x - 1)(-1)^n.
        \end{equation} 

    In this note, we consider the $n$-cycle graph to be the planar graph associated with the $(2,n)$-torus knot diagram (\Figs{fig:link}~\subref{subfig:ncyclegraph}). 
    
        \begin{table}[!ht]
            \centering
            $\begin{array}{c|rrrrrrrrrrrrrrrrrrrrr}
            n\ \backslash\ k & 0 & 1 & 2 & 3 & 4 & 5 & 6 & 7\\
            \hline
            1 & 0\\  
            2 & 0 & -1 &   1\\
            3 & 0 &  2 &  -3 &   1\\
            4 & 0 & -3 &   6 &  -4 &   1\\
            5 & 0 &  4 & -10 &  10 &  -5 &  1\\
            6 & 0 & -5 &  15 & -20 &  15 & -6 &  1\\
            7 & 0 &  6 & -21 &  35 & -35 & 21 & -7 &  1
            \end{array}$
            \caption{First $7$ rows in \seqnum{A137396}.}
            \label{tab:A137396}
        \end{table}
    
    \item Row $n$ in \seqnum{A300453} is generated by 
    
        \begin{equation}\label{eq:gen}
            xK_n(1,1,x) = (x + 1)^n + x^2 - 1.
        \end{equation}
    
    We referred to the polynomial in \eqref{eq:gen} as ``{generating polynomial}'' \cite{Ramaharo}. This is, in fact, the expression of the bracket evaluated at the shadow diagram (\Figs{fig:link}~\subref{subfig:torusshadow}).

        \begin{table}[!ht]
            \centering
            $\begin{array}{c|rrrrrrrrrrrrr}
        	n\ \backslash\ k	&0	&1	&2	&3	&4	&5	&6	&7\\
        	\midrule
        	0	&0	&0	&1	&	&	&	&	&\\
        	1	&0	&1	&1	&	&	&	&	&\\
        	2	&0	&2	&2	&	&	&	&	&\\
        	3	&0	&3	&4	&1	&	&	&	&\\
        	4	&0	&4	&7	&4	&1	&	&	&\\
        	5	&0	&5	&11	&10	&5	&1	&	&\\
        	6	&0	&6	&16	&20	&15	&6	&1	&\\
        	7	&0	&7	&22	&35	&35	&21	&7	&1
        	\end{array}$
            \caption{First $8$ rows in \seqnum{A300453}.}
            \label{tab:A300453}
        \end{table}
\end{itemize}

We show in the next section the connection between these polynomials.

\section{Construction and interpretation}\label{sec:interpretation}
\subsection{Bracket polynomial}
The bracket polynomial for the knot diagram $K$ is defined by
    \begin{equation}\label{eq:bpdfn}
    	\langle K \rangle = K(A,B,d) = \sum_{s}^{}\langle K|s \rangle d^{|s| - 1},
    \end{equation}
where $ \langle K|s \rangle $ denotes the product of the splitting variables ($ A $ and $ B $) associated with the state $ s $, and $|s|$ denotes the number of circles (or loops) in $ s $.

\begin{figure}[!ht]
	\centering
	\includegraphics[width=0.3\linewidth]{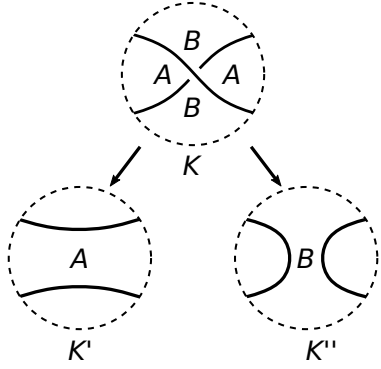}
	\caption{$ A $-split and $ B $-split.}
	\label{fig:statesab}
\end{figure}

Formula \eqref{eq:bpdfn} can also be expressed as
\begin{itemize}
    \item $\langle K  \rangle = A\langle K' \rangle + B\langle K'' \rangle$,
    \item $\langle \bigcirc \bigcirc\cdots\bigcirc\rangle = d^{k - 1}$ (disjoint union of $k$ circles), 
\end{itemize}
where  $K' $ and $ K'' $ are obtained from $ K $ by performing $ A $ and $ B $ splits at a given crossing in $ K $, see \Figs{fig:statesab}. Last formula reads as well as $\langle \bigcirc\ K\rangle = d\langle K\rangle$ and $\langle \bigcirc\rangle = 1$.

For example,
$ \langle T_n \rangle = (A + Bd)^n$, where $ T_n := \includegraphics[width=0.2\linewidth,valign=c]{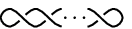} $ (with $ n $ half-twists).

Indeed, we have 
    \begin{align*}
        \left< \includegraphics[width=0.2\linewidth,valign=c]{twist} \right> & = A\left< \includegraphics[width=0.2\linewidth,valign=c]{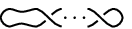}  \right> + B\left< \includegraphics[width=0.2\linewidth,valign=c]{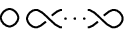} \right>\\
        & = A\langle T_{n-1}\rangle + B\langle\bigcirc\ T_{n-1}\rangle `\\
        & = \left(A + Bd\right) \langle T_{n-1}\rangle\\
        & = \cdots\\
        & = (A + Bd)^n.
    \end{align*}

Now, formula \eqref{eq:brack} is straightforward if we notice that,
    \begin{equation}\label{key}
        \left< K_n\right> = A \left< K_{n-1}\right> + B\left< T_{n-1}\right>, \ \textit{with}\  \left< K_0\right> = d.
    \end{equation}
Finally, set $ B = A^{-1} $ and $ d = -A^{-2} - A^{2} $ so that the bracket is invariant under Reidemeister moves II and III \cite[p.\ 31--33]{Kauffman}. For $n = 0, 1, \ldots, 5$ the corresponding bracket polynomial reads
    \begin{align*}
        \left< K_0\right> & = -A^2 - \dfrac{1}{A^2}\\
        \left< K_1\right> & = -A^3 \\
        \left< K_2\right> & = -A^4 - \dfrac{1}{A^4}\\
        \left< K_3\right> & = -A^5 - \dfrac{1}{A^3} + \dfrac{1}{A^7}\\
        \left< K_4\right> & = -A^6 - \dfrac{1}{A^2} + \dfrac{1}{A^6} - \dfrac{1}{A^{10}}\\
        \left< K_5\right> & = -A^7 - \dfrac{1}{A}   + \dfrac{1}{A^5} - \dfrac{1}{A^9}  + \dfrac{1}{A^{13}}.
    \end{align*}
The $n$-th row polynomial of triangle in \seqnum{A123192} is then obtained by multiplying the bracket $\langle K_n \rangle$ by  $ A^{|3n - 2|} $.

\subsection{Chromatic polynomial}
\begin{figure}[!ht]
	\centering
	\includegraphics[width=0.3\linewidth]{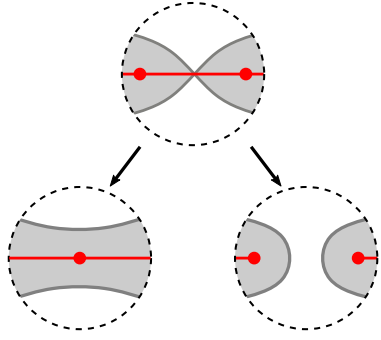}
	\caption{Edge contraction and deletion.}
	\label{fig:chromaticstate}
\end{figure}

In the present framework, the splits of type $ A $ and $ B $ may be regarded in terms of graph as the ``contraction'' and ``deletion'' operations, respectively, as shown in \Figs{fig:chromaticstate}. Kaufman refers to the resulting states as ``chromatic states'' \cite[p.\ 353--358]{Kauffman}.

Let $n\geq 1$, and let $ G(K_n,x) = xK(A,B,x) = \sum_{s}^{}\langle K|s \rangle x^{|s|}$. By \eqref{eq:chrom}, we can rewrite \eqref{eq:bpdfn} as
    \begin{equation}
    	G(K_n,x) = \sum_{s}^{}(-1)^{A(s)}\left(x^\frac{1}{2}\right)^{B(s)}\left(x^\frac{1}{2}\right)^{|s|}, 
    \end{equation}
where $A(s)$ and $B(s)=n - A(s)$ are the number of $ A $-splits and $ B $-splits in the state $s$, respectively. Furthermore, we have  $|s| = 2$ if $B(s) = 0$, and $|s| = B(s)$ otherwise \cite{Ramaharo}.  Hence
\begin{itemize}
    \item $G(K_1,x) = (-1)^{1}\left(x^\frac{1}{2}\right)^{0}\left(x^\frac{1}{2}\right)^{2} + (-1)^{0}\left(x^\frac{1}{2}\right)^{1}\left(x^\frac{1}{2}\right)^{1} = 0$;
    \item and for $n\geq2$,
\end{itemize}
    \begin{align*}
    	G(K_n,x) & = \sum_{s}^{}(-1)^{n - B(s)}x^{\frac{1}{2}{\left(B(s) + |s|\right)}}\\
                 & = (-1)^{n}x + \sum_{|s|\geq1} (-1)^{n - B(s)}x^{B(s)}\\
                 & = \sum_{s} (-1)^{i(s)}x^{B(s)},
    \end{align*}

where $ i(s)$ is the number of ``interior vertices'' in the chromatic state $ s $ \cite[p.\ 358]{Kauffman} and 
 $ B(s) $ matches the number of shaded components in $ s $ \cite{Ramaharo} (with $ i(s) = 1 $ if $ |s| = 1 $, and $ i(s) = n - B(s) $ otherwise).

\subsection{Generating polynomial}
Now, what if we evaluate the bracket polynomial at the shadow diagram? Recall that a shadow is a knot diagram without under or over-crossing information. In such case, it is natural to set $ A = B = 1 $ in \eqref{eq:bptn}. The bracket becomes
    \begin{equation*}
        K_n(1,1,d) = \dfrac{(d + 1)^n + d^2 - 1}{d}, 
    \end{equation*}
and formula \eqref{eq:gen} implies
    \begin{equation*}
        xK_n(1,1,x) = \sum_{s} x^{|s|} =\sum_{k} s(n,k) x^k,
    \end{equation*}
where $s(n,k)$ is the number of states having exactly $ k $ circles  \cite{Ramaharo}.

\bigskip
\hrule
\bigskip

\noindent 2010 {\it Mathematics Subject Classification}: 05A10; 57M25.

\bigskip
\hrule
\bigskip

\noindent (Concerned with sequences
\seqnum{A123192}, \seqnum{A137396} and \seqnum{A300453}.)
\bigskip
\hrule
\bigskip

\end{document}